# ANALYTIC NON-LINEARIZABLE UNIQUELY ERGODIC DIFFEOMORPHISMS ON $\mathbb{T}^2$

Maria Saprykina

ABSTRACT. In this paper we study the behavior of diffeomorphisms, contained in the closure $\overline{\mathcal{A}_\alpha}$ (in the inductive limit topology) of the set $\mathcal{A}_\alpha$ of real-analytic diffeomorphisms of the torus $\mathbb{T}^2$, conjugated to the rotation $R_\alpha : (x, y) \mapsto (x + \alpha, y)$ by an analytic measure-preserving transformation. We show that for a generic $\alpha \in [0, 1]$, $\overline{\mathcal{A}_\alpha}$ contains a dense set of uniquely ergodic diffeomorphisms. We also prove that $\overline{\mathcal{A}_\alpha}$ contains a dense set of diffeomorphisms that are minimal and non-ergodic.

## Introduction

Consider the set Diff($\mathbb{T}^2$) of analytic diffeomorphisms of the torus homotopic to the identity (where $\mathbb{T}^2 = \mathbb{R}^2/\mathbb{Z}^2$). We provide Diff($\mathbb{T}^2$) with the inductive limit topology, induced by the supremum norms of analytic functions over the complex neighborhoods of $\mathbb{T}^2$ (see Section 1.1).

This paper is devoted to the study of the following subsets of Diff($\mathbb{T}^2$): for any $\alpha \in [0, 1]$, $\mathcal{A}_\alpha$ is the set of analytic diffeomorphisms $F : \mathbb{T}^2 \to \mathbb{T}^2$, which are analytically conjugated to the rotation $R_\alpha : (x, y) \mapsto (x + \alpha, y)$, i.e. such that there exists an analytic area preserving diffeomorphism $T$ of the torus such that $F = T^{-1} \circ R_\alpha \circ T$.

For any $\alpha$, a generic small exact symplectic perturbation of a diffeomorphism from $\mathcal{A}_\alpha$ exhibits the whole spectrum of models of behavior: there is a set of large measure foliated by invariant circles with uniquely ergodic motion on each of them (this is described by the KAM-theorem), hyperbolic and elliptic periodic points, elliptic islands.

Here we investigate a special type of such perturbations, namely those contained in the closure of $\mathcal{A}_\alpha$. The set $\overline{\mathcal{A}_\alpha}$ of such perturbations is "small" in Diff($\mathbb{T}^2$). Indeed, all the diffeomorphisms $F \in \overline{\mathcal{A}_\alpha}$ are formally conjugated to the rotation in the following sense: there exists a formal Fourier series $T$ such that $T^{-1} \circ F \circ T = R_\alpha$. This is due to the fact that each of the Fourier coefficients of the normalizing transformation $T$ is a polynomial, depending only on finitely many Fourier coefficients of $F$. At the same time, the formal normal form of an arbitrary exact symplectic perturbation of a diffeomorphism from $\mathcal{A}_\alpha$ is $(x, y) \mapsto (x + \alpha + f(y), y)$. This remark indicates that $\overline{\mathcal{A}_\alpha}$ has infinite codimension in Diff($\mathbb{T}^2$).

Opposite to the case of generic perturbations in Diff($\mathbb{T}^2$), it is probable that the properties of the sets $\overline{\mathcal{A}_\alpha}$ are different, depending on the arithmetical properties of $\alpha$. For a Diophantine (or Bruno) $\alpha$, one can believe that the analogy with the following result of Rüssmann [Rüs], [B], [E1], holds true. It claims, in particular, that if a real analytic Hamiltonian transformation in a neighborhood of an elliptic fixed point is

Typeset by $\mathcal{A}\mathcal{M}\mathcal{S}$-TEX

formally conjugated to the *linear* normal form, then the normalizing transformation actually converges, provided that the eigenvalues satisfy a Diophantine (or Bruno) condition. In parallel to this result, it seems likely that, in the case of a Diophantine $\alpha$, for any $F \in \mathcal{A}_\alpha$ there exists a neighborhood $U(F)$ such that $U(F) \bigcap \overline{\mathcal{A}_\alpha} \subset \mathcal{A}_\alpha$.

On the other hand, we show that for a generic set of real numbers $\alpha$ (by "generic" we mean "a dense $G_\delta$-set"), namely for those $\alpha$ that are "well enough approximated by rational numbers", there is a dense set of non-linearizable diffeomorphisms in $\overline{\mathcal{A}_\alpha}$, i.e. such that there is no homeomorphism of the torus, conjugating it to a rotation. But much more is indeed true. Let us recall that a transformation of a compact Hausdorff space is called uniquely ergodic if it leaves invariant a unique probability measure on this space.

**Theorem A.** *For a generic set of real numbers $\alpha$, $\overline{\mathcal{A}_\alpha}$ contains a dense set of uniquely ergodic diffeomorphisms.*

Speaking of uniquely ergodic diffeomorphisms, one expects their genericity; but since the set of real-analytic diffeomorphisms with our topology is not a Baire space, genericity is not stronger than density in it. In the last section we discuss the question of genericity of unique ergodicity in particular subspaces of $\overline{\mathcal{A}_\alpha}$.

Unique ergodicity with respect to the Lebesgue measure implies minimality, hence it follows from Theorem A that for a generic set of $\alpha$, $\overline{A_\alpha}$ contains a dense set of minimal diffeomorphisms. In addition to this, we prove the following result.

**Theorem B.** *For a generic set of real numbers $\alpha$, $\overline{\mathcal{A}_\alpha}$ contains a dense set of diffeomorphisms that are minimal and non-ergodic.*

As a corollary from Theorem A we obtain the fact that $\overline{A_\alpha}$ contains a dense set of diffeomorphisms with zero topological entropy. The fact, that on a two-dimensional manifold minimality implies vanishing of the topological entropy, was proven by A. Katok in [K].

Our constructions are specifically adapted to the case of the torus. An interesting related question is: What is the typical behavior of a real analytic measure preserving diffeomorphism near its closed periodic curve? Near an elliptic fixed point?

One of the motivations for the present work is the following. Anosov and Katok in their well-known paper [AK] of 1970 presented a set of examples of non-linearizable $C^\infty$-diffeomorphisms enjoying ergodic, weak mixing, and other statistical properties. The construction can be performed on an arbitrary smooth manifold $M$, supporting a periodic flow. Moreover, it is shown that in the closure of the set of diffeomorphisms of $M$, contained in smooth, measure preserving periodic flows on $M$, both ergodic and weakly mixing diffeomorphisms are generic. Their construction is based on the method of fast cyclic approximations by periodic diffeomorphisms, that relies on the existence of functions with compact support in $C^\infty$. Similar methods permit A.Fathi and M.Herman [FH] to prove the existence of smooth minimal and uniquely ergodic diffeomorphisms on smooth manifolds, supporting a periodic flow.

Work in the analytic category demands completely different techniques, and the corresponding set of examples does not yet exist in the whole generality. In this paper we present a part of the analogous treatment of the analytic case.

Let us briefly recall some of the related results in the analytic category. There is an example, due to H.Furstenberg [M], [KH], of an analytic area-preserving dif-feomorphism of $\mathbb{T}^2$, that is minimal, but not ergodic. Our method allows to obtain this result as well (Theorem 1.2) (moreover, in Theorem B we present a dense set of such diffeomorphisms).

Uniquely ergodic examples constructed in the present work have a form of a skew-product with an irrational rotation in the base: $G(x, y) = (x + \alpha, y + g(x, y))$. In the work [Fu], H.Furstenberg studies this type of analytic skew-products, and formulates a criterion of unique ergodicity for them. He also shows that for such diffeomorphisms ergodicity in the usual sense implies unique ergodicity. We give a constructive procedure for defining uniquely ergodic skew-products.

Ergodic properties of analytic *skew systems* were studied by L.H.Eliasson in [E2]. Namely, on $\mathbb{T}^d \times SO(3, \mathbb{R})$ he considers the systems of the form

$$\begin{cases} \dot{X} = F(x)X \\ \dot{x} = \omega, \end{cases}$$

where $F : \mathbb{T}^d \to o(3, \mathbb{R})$ is real analytic, and $\omega \in \mathbb{R}^d$ satisfies a Diophantine condition. L.H.Eliasson proves that, for a generic $F$ close to constant coefficients, this system is uniquely ergodic.

Weak mixing property for reparametrized irrational translation flows on the torus $\mathbb{T}^n$, $n \geq 2$, was investigated in its full generality by B.Fayad, [F1]. In [F2], he presents a study of mixing property for reparametrized irrational flows on the torus $\mathbb{T}^3$ (or any $\mathbb{T}^n$, $n \geq 3$). Mixing examples are obtained using analytic time changes in the irrational translation flow.

1. DEFINITIONS AND PLAN OF THE PROOF

**1.1. Topology on $\mathcal{A}_\alpha$ and $\mathcal{A}_\alpha^r$.** The set $\mathcal{A}_\alpha$ under consideration lies in the set $\text{Diff}(\mathbb{T}^2)$ of analytic diffeomorphisms of the torus, homotopic to the identity. Let us introduce the topology $\tau$ on the latter.

For any $r > 0$ we define $C_r^\omega(\mathbb{T}^2)$ as the set real analytic functions on $\mathbb{R}^2$, $\mathbb{Z}^2$-periodic, that can be extended to holomorphic functions on $A^r = \{|\text{Im}\, x|, |\text{Im}\, y| < r\}$. On this space we use the uniform norm $|f|_r = \sup_{A^r} |f(x, y)|$.

Elements of $\text{Diff}(\mathbb{T}^2)$ are homotopic to the identity, and hence have a lift of type $F(x, y) = (x + f_1(x, y), y + f_2(x, y))$ with $f_i$ analytic on $\mathbb{R}^2$ and $\mathbb{Z}^2$-periodic. Consider the subspace $\text{Diff}_r(\mathbb{T}^2)$ of $\text{Diff}(\mathbb{T}^2)$, consisting of those diffeomorphisms, for whose lift it holds: $f_i \in C_r^\omega(\mathbb{T}^2)$, $i = 1, 2$. There is a natural isomorphism between $\text{Diff}_r(\mathbb{T}^2)$ and the Banach space $\mathcal{D}_r$ of pairs of functions in $C_r^\omega(\mathbb{T}^2)$, with the topology defined by the supremum norms. We endow $\text{Diff}_r(\mathbb{T}^2)$ with the topology $\tau_r$, brought from $\mathcal{D}_r$ by this isomorphism.

For diffeomorphisms $F$, $G$ in $\text{Diff}_r(\mathbb{T}^2)$ we shall use the distance, generating the same topology:
$$|F - G|_r = \max_{i=1,2}\{|f_i - g_i|_r\}.$$

Space $\text{Diff}(\mathbb{T}^2)$ is isomorphic to the union $\mathcal{D} = \bigcup \mathcal{D}_r$. Here $(\mathcal{D}_r)_r$ is a growing (when $r$ goes to zero) family of Banach spaces. For every pair $r > s$ we have a continuous linear injection $i_{rs} : \mathcal{D}_r \to \mathcal{D}_s$, the image $i_{rs}(\mathcal{D}_r)$ is dense in $\mathcal{D}_s$, and the map $i_{rs}$ is compact (i.e., the images of balls in $\mathcal{D}_r$, $i_{rs}B_r(f, R)$, are pre-compact in $\mathcal{D}_s$). In this case we can endow $\mathcal{D}$ with the inductive limit topology: a set $O \subset \mathcal{D}$ is open in $\mathcal{D}$ if and only if for any $r > 0$ the set $O \cap \mathcal{D}_r$ is open in $\mathcal{D}_r$. A good description of this topology can be found in [L] (Appendix 2). Alternatively, one can define the inductive limit topology $\tau$ in the usual way (see, for example, [R]), and prove the equivalence with the above definition in our case. Finally, we define the topology $\tau$ on Diff($\mathbb{T}^2$) bringing the above inductive limit topology to Diff($\mathbb{T}^2$) by the natural isomorphism.

We would like to stress that $\tau$ does not make Diff($\mathbb{T}^2$) into a Baire space, and we cannot, unfortunately, speak about genericity in it.

The set $\mathcal{A}_\alpha$ gets the subspace topology, generated by $\tau$.

Up to the last section we shall study the following subset $\mathcal{A}_\alpha^r$ of $\mathcal{A}_\alpha \cap \text{Diff}_r(\mathbb{T}^2)$ (for an arbitrary fixed $r > 0$):

$$\mathcal{A}_\alpha^r = \{F \in \mathcal{A}_\alpha \cap \text{Diff}_r(\mathbb{T}^2) \mid F = T^{-1} \circ R_\alpha \circ T \text{ for some real analytic,}$$
$$\text{area preserving } T, \text{ such that } T^{-1} \text{ is analytic in a neighborhood of } \overline{R_\alpha \circ T(A^r)}\}$$

with the topology induced by $\tau_r$. $\mathcal{A}_\alpha^r$ is a Banach space, and we shall prove genericity of uniquely ergodic diffeomorphisms and density of minimal non-ergodic ones in $\overline{\mathcal{A}_\alpha^r}$. Note that $\mathcal{A}_\alpha^r \subsetneq \mathcal{A}_\alpha \cap \text{Diff}_r(\mathbb{T}^2)$, and we do not know how to prove the same statement for this bigger set.

The set $\mathcal{A}_\alpha$ equals $\bigcup_{r>0} \mathcal{A}_\alpha^r$. By the construction of our topologies, a sequence, converging in $\mathcal{A}_\alpha^r$, converges to the same limit in $\mathcal{A}_\alpha$. Hence, $\overline{\mathcal{A}_\alpha^r} \subset \overline{\mathcal{A}_\alpha}$ for any $r$ (closures are taken in the corresponding topologies), and density of both minimal non-ergodic and uniquely ergodic diffeomorphisms in $\overline{\mathcal{A}_\alpha^r}$ implies their density in $\overline{\mathcal{A}_\alpha}$ with topology induced by $\tau$.

In what follows we shall not distinguish between a diffeomorphism of the torus and its lifts when this does not lead to a confusion.

**1.2. Main construction.** The method is based on the following idea. In Lemma 2.1 we show that, provided that a real number $\alpha$ is "well enough approximated by rationals", for any $\varepsilon > 0$ and $r > 0$, one can find an area preserving real-analytic transformation $T$, homotopic to the identity, such that the dynamics of $G := T^{-1} \circ R_\alpha \circ T$ is very different from that of $R_\alpha$ (in particular, the invariant curves are very far from circles), while the difference $|G - R_\alpha|_r$ is less than $\varepsilon$.

In the process of the construction we shall step by step produce a converging sequence of analytic diffeomorphisms $\{G_n\}$, each time conjugating $R_\alpha$ by "wilder and wilder" transformations $T_n$: $G_n = T_n^{-1} \circ R_\alpha \circ T_n$ so that $|G_n - G_{n-1}|_r$ go to zero very fast. The desired ergodic diffeomorphism is defined as a limit

$$G = \lim_{n \to \infty} G_n.$$

To make this idea work, the value of $\alpha$ has to be chosen close to certain rational numbers $\frac{p_n}{q_n}$, $n \in \mathbb{N}$. (We always assume that $(p_n, q_n) = 1$, i.e. $p_n$ and $q_n$ are relatively prime.) We have chosen a way of presentation where we do not fix the value of $\alpha$ in advance, but construct it as a limit of the inductive procedure. Approximating diffeomorphisms $G_n$ will thus depend on a real parameter $\alpha$. We express this dependence as $G_n(\alpha; x, y)$ or $G_n(\alpha)$. At the $n$-th step we shall construct $G_n(\alpha)$ and choose a closed interval $I_n$, centered at $\frac{p_n}{q_n}$, $I_n \subset I_{n-1}$ such that the desired estimates for $G_n(\alpha)$ (in particular, $|G_n(\alpha) - G_{n-1}(\alpha)|_r$ small) hold for any $\alpha \in I_n$. The number $\hat{\alpha}$, providing all the necessary estimates, is $\hat{\alpha} = \bigcap_n I_n = \lim_n \frac{p_n}{q_n}$. If $q_n$ grow with $n$, this number is irrational.

It is important to note that we do not need any lower bounds on $|I_n|$-s, and we do not search for the sharp conditions on the parameters. Numbers $\hat{\alpha}$ that can be obtained as limits of the above construction, are characterized only by "good approximation by rationals", therefore (see Lemma 7.2), such numbers $\hat{\alpha}$ form a $G_\delta$-set in $[0, 1]$.

The only parameters of the construction are a sequence of nested closed intervals $I_n \subset I_{n-1}$, centered at rational points $\frac{p_n}{q_n}$, $n \in \mathbb{N}$ with the corresponding lengths $|I_n|$, and a sequence of positive numbers $c_n$. The numbers $p_n$, $q_n$, $c_n$ and $|I_n| \in \mathbb{R}$ are chosen at the $n$-th step of the construction.

Define $G_0 = R_\alpha$ and $I_0 = [0, 1]$, and denote $r_n = 10^n$, $\varepsilon_n = \frac{1}{10^{n+1}}$. Let us describe the $n$-th step of the construction more precisely (still, in the full details it appears only in the proof). Suppose that we have a nested sequence of closed intervals $I_j$, $j = 1, \ldots n-1$, centered at rational points $\frac{p_j}{q_j}$, a rational point $\frac{p_n}{q_n} \in I_{n-1}$, and a sequence of positive real numbers $c_j$, $j = 1, \ldots n$. Let $\Phi_n$ be the flow map over the time $\frac{1}{2q_n}$ of the Hamiltonian vector field with the Hamiltonian function

$$H_n = y + \frac{c_n}{2} \sin 2\pi q_n x. \tag{1.1}$$

Explicitly,

$$\Phi_n : (x, y) \mapsto \left(x + \frac{1}{2q_n},\ y + c_n \sin(2\pi q_n x)\right). \tag{1.2}$$

Note that $\Phi_n$ is area preserving, since it is a Hamiltonian flow map. Let $\Phi_j$ be defined for all $j \leq n$ by the above formula with $j$ instead of $n$. The conjugating transformation at the $n$-th step has the form

$$T_n = \Phi_n \circ \Phi_{n-1} \circ \ldots \circ \Phi_1,$$

and the $n$-th approximation $G_n(\alpha)$ is defined as the composition

$$G_n(\alpha) := T_n^{-1} \circ R_\alpha \circ T_n. \tag{1.3}$$

We shall choose $I_n \subset I_{n-1}$ centered at $\frac{p_n}{q_n}$ so small that for any $\alpha \in I_n$

$$|G_n(\alpha) - G_{n-1}(\alpha)|_{r_n} < \varepsilon_n. \tag{1.4}$$

The possibility of this is proven in Proposition 2.1. We choose $\frac{p_{n+1}}{q_{n+1}} \in I_n$ with a large denominator, and repeat the procedure iteratively. Then for $\hat{\alpha} = \bigcap_n I_n$ the analytic limit $G(\hat{\alpha}) = \lim G_n(\hat{\alpha})$ exists. Moreover,

$$|G(\hat{\alpha}) - R_{\hat{\alpha}}|_{10} \leq |G_1(\hat{\alpha}) - R_{\hat{\alpha}}|_{r_1} + \sum_{n=2}^{\infty} |G_n(\hat{\alpha}) - G_{n-1}(\hat{\alpha})|_{r_n} < \frac{1}{10}.$$

We shall show that under certain conditions on the decay of $|I_n|$ and growth of $c_n$, the limit $G(\hat{\alpha})$ is analytic non-linearizable (minimal non-ergodic, uniquely ergodic).

**1.3. Modification used to obtain density.** The above argument permits us, given $r$ and $\varepsilon > 0$, to find a number $\hat\alpha$ and a diffeomorphism $G(\hat\alpha) \in \overline{\mathcal{A}^r_{\hat\alpha}}$ with the required properties, such that $|G(\hat\alpha) - R_{\hat\alpha}|_r < \varepsilon$. In order to get a dense set of such examples in $\overline{\mathcal{A}^r_{\hat\alpha}}$ for an appropriate $\hat\alpha$, we modify the above construction a little. Namely, together with $\hat\alpha$, we construct not one, but an infinite sequence of diffeomorphisms $G_{(m)}(\hat\alpha)$ with desired properties, pairwise conjugated by area-preserving analytic transformations (having, therefore, the same ergodic properties), such that

$$|G_{(m)}(\hat\alpha) - R_{\hat\alpha}|_{10^m} < \frac{1}{10^m}.$$

In the last section we show that this is enough to imply density in $\overline{\mathcal{A}^r_{\hat\alpha}}$.

Diffeomorphisms $G_{(m)}(\hat\alpha)$ are obtained in the same way as $G(\hat\alpha)$, but with $c_1 = \cdots = c_m = 0$. In other words, for any $n$ we consider

$$T_{m,n} = \Phi_n \circ \Phi_{n-1} \circ \ldots \circ \Phi_m \quad \text{for } m \leq n, \quad \text{and } T_{n+1,n} = \operatorname{Id},$$

and

(1.5) $\qquad G_{m,n}(\alpha) = T^{-1}_{m,n} \circ R_\alpha \circ T_{m,n} \quad \text{for } m \leq n, \quad \text{and } G_{n+1,n}(\alpha) = R_\alpha.$

If the limit of $G_{m,n}(\alpha)$ exists, denote for any fixed $m$ ($n = m, \ldots$)

$$\lim_{n\to\infty} G_{m,n}(\alpha) = G_{(m)}(\alpha).$$

In particular, $T_{1,n} = T_n$, $G_{1,n}(\alpha) = G_n(\alpha)$. For $G_{(1)}(\alpha)$ we keep the notation $G(\alpha)$. Now at the $n$-th step of the construction we shall choose $I_n$ so small that, in addition to (1.4), for any $m \leq n$ we have:

$$|G_{m,n}(\alpha) - G_{m,n-1}(\alpha)|_{r_n} < \varepsilon_n$$

for all $\alpha \in I_n$. The possibility of this choice follows from Proposition 2.1. Then for $\hat\alpha = \bigcap_n I_n$ we have, for any $m$, a sequence of diffeomorphisms $G_{m,n}(\alpha)$, converging to a non-linearizable (minimal, ergodic) analytic limit $G_{(m)}$, and

$$|G_{(m)}(\hat\alpha) - R_{\hat\alpha}|_{10^m} \leq \sum_{n=m}^{\infty} |G_{m,n}(\hat\alpha) - G_{m,n-1}(\hat\alpha)|_{r_n} < \frac{1}{10^m}.$$

*Remark 1.1.* Clearly, for any $n$ and $m \leq n$, $G_{m,n}(\alpha)$ is conjugated to $G_n(\alpha)$: $G_{m,n}(\alpha) = S_m^{-1} \circ G_n(\alpha) \circ S_m$, where $S_m^{-1} = \Phi_{m-1} \circ \ldots \circ \Phi_1$. Therefore, for a given $\alpha$ either for each $m$ the sequence $G_{m,n}(\alpha)$ converges with $n \to \infty$, or it does not converge for any $m$. If for some $\alpha$ the sequences do converge, then

$$G_{(m)}(\alpha) = S_m^{-1} \circ G(\alpha) \circ S_m.$$

Since $S_m$ is analytic and area preserving, the limits $G_{(m)}(\alpha)$ have the same ergodic properties for all $m$: they are either minimal (respectively, ergodic or uniquely ergodic) for all $m$, or not.

**1.4. Plan of the article.** Section 2 is dedicated to the preliminary technical work. In Section 3, with the help of the above construction, we produce non-linearizable analytic diffeomorphism arbitrary close to the rotation. "Non-linearizable" means here that there is no homeomorphism of $\mathbb{T}^2$ such that $G(\alpha) = T^{-1} \circ R_\alpha \circ T$.

**Theorem 1.1.** *Suppose that $\sum_n c_n = \infty$, $\alpha = \lim \frac{p_n}{q_n}$. If $|\frac{p_n}{q_n} - \alpha|$ decay sufficiently fast, then for any $r$ and $\varepsilon > 0$ there exists a non-linearizable analytic diffeomorphism $G_{r,\varepsilon} \in \overline{\mathcal{A}_\alpha^r}$ such that $|G_{r,\varepsilon} - R_\alpha|_r < \varepsilon$.*

In Section 4 we chose the parameters to obtain minimal non-ergodic examples.

**Theorem 1.2.** *Let $c_n = \frac{1}{n}$ for all $n \in \mathbb{N}$, and $\alpha = \lim \frac{p_n}{q_n}$. If $|\frac{p_n}{q_n} - \alpha|$ decay sufficiently fast, then for any $r$ and $\varepsilon > 0$ there exists a minimal non-ergodic analytic diffeomorphism $G_{r,\varepsilon} \in \overline{\mathcal{A}_\alpha^r}$ such that $|G_{r,\varepsilon} - R_\alpha|_r < \varepsilon$.*

Sections 5 and 6 contain the proof of the following result.

**Theorem 1.3.** *Suppose that $c_n$ grow sufficiently fast, and $\alpha = \lim \frac{p_n}{q_n}$. If $|\frac{p_n}{q_n} - \alpha|$ decay sufficiently fast, then for any $r$ and $\varepsilon > 0$ there exists a uniquely ergodic analytic diffeomorphism $G_{r,\varepsilon} \in \overline{\mathcal{A}_\alpha^r}$ such that $|G_{r,\varepsilon} - R_\alpha|_r < \varepsilon$.*

The last section concerns density of minimal non-ergodic and uniquely ergodic diffeomorphisms in $\overline{\mathcal{A}_\alpha}$ for an appropriate $\alpha$. We show as well that these "appropriate" $\alpha$ are generic in $[0,1]$. This will finish the proof of Theorems A and B.

At the end of the last section we discuss genericity of uniquely ergodic diffeomorphisms in $\overline{\mathcal{A}_\alpha^r}$ for an appropriate $\alpha$.

## 2. Important proposition

In this section we show that the choice of a sufficiently small interval $I_n$ (centered at $\frac{p_n}{q_n}$) implies that the difference $G_{n-1}(\alpha) - G_n(\alpha) = \hat{G}_n(\alpha)$ becomes small for all $\alpha \in I_n$. This statement follows from the continuity of the conjugation. Indeed, by (1.2) and (1.3), for $\alpha = \frac{p_n}{q_n}$ the difference $G_{n-1}(\alpha) - G_n(\alpha)$ equals zero. One would expect that $G_n(\alpha)$ is still close to $G_{n-1}(\alpha)$ for $\alpha$ sufficiently close to $\frac{p_n}{q_n}$.

**Proposition 2.1.** *Suppose that $T \in \bigcap_{r>0} \text{Diff}_r(\mathbb{T}^2)$, and let*

$$G(\alpha) = T^{-1} \circ R_\alpha \circ T.$$

*Then for any $\frac{p_n}{q_n} \in [0,1]$, $c_n$, $r \geq 0$ and $\varepsilon > 0$ there exists an interval $I = I(T, p_n, q_n, c_n, r, \varepsilon)$, centered at $\frac{p_n}{q_n}$, such that for any $\alpha \in I$ the diffeomorphism*

(2.1) $$G_n(\alpha) = T_n^{-1} \circ R_\alpha \circ T_n, \quad T_n = \Phi_n \circ T$$

*with $\Phi_n$ defined by (1.2), satisfies*

$$|G_n(\alpha) - G(\alpha)|_r < \varepsilon.$$

**Corollary 2.1.** *Suppose that $G(\alpha)$ and $T$ are like in Proposition 2.1. Then for any $\frac{p_n}{q_n} \in [0,1]$, $c_n$, $r \geq 0$, $\varepsilon > 0$ and $\tau \in \mathbb{N}$ there exists an interval $I = I(T, p_n, q_n, c_n, r, \varepsilon, \tau)$, centered at $\frac{p_n}{q_n}$, such that for any $\alpha \in I$ the diffeomorphism $G_n(\alpha)$, defined by (2.1), satisfies*

$$\max_{i=0,\ldots\tau} |G^i(\alpha) - G_n^i(\alpha)|_r < \varepsilon. \tag{2.2}$$

The following lemma constitutes the kernel of the proof of Proposition 2.1.

**Lemma 2.1.** *Let $\Phi_n$ be defined by (1.2). Then for any $\frac{p_n}{q_n} \in [0,1]$, $c_n$, $r \geq 0$, and $\varepsilon > 0$ there exists an interval $I = I(p_n, q_n, c_n, r, \varepsilon)$, centered at $\frac{p_n}{q_n}$, such that for any $\alpha \in I$ the difference*

$$\Phi_n^{-1} \circ R_\alpha \circ \Phi_n - R_\alpha =: \hat{F}_n(\alpha) \tag{2.3}$$

*is small in $r$-metric: $|\hat{F}_n(\alpha)|_r < \varepsilon$.*

*Proof.* Let us first estimate the difference $\tilde{F}_n(\alpha) := (R_\alpha \circ \Phi_n - \Phi_n \circ R_\alpha)$ (seen as a complex-analytic function with $\mathbb{Z}^2$-periodic real part). For this, let us set $N = \max_{\alpha \in [0,1]} |\tilde{F}_n(\alpha)|_r$, and for $r_1 = (\max_{A^r} |\mathrm{Im}(\Phi_n)| + N)$, let $M = |D\Phi_n^{-1}|_{r_1}$. We show first, that for a sufficiently small interval $I$ around $\frac{p_n}{q_n}$, $\max_{\alpha \in I} |\tilde{F}_n(\alpha)|_r$ is less than $\frac{\varepsilon}{M}$.

Write down $\tilde{F}_n(\alpha)$ explicitly, omitting the indices (use (1.2)):

$$|\tilde{F}(\alpha)|_r = \left| \left( \frac{1}{2q} + \alpha, c\sin(2\pi q x) \right) - \left( \frac{1}{2q} + \alpha, c\sin(2\pi q(x+\alpha)) \right) \right|_r$$
$$= c |\sin(2\pi q x) - \sin(2\pi q(x+\alpha))|_r \leq cC|1 - e^{2\pi i q \alpha}|,$$

where $C = e^{2\pi q r}$. Now it is enough to take the length of $I$ such that $0 < |I| < \frac{\varepsilon}{4\pi MCqc}$ and $I \subset [0,1]$. Indeed, any real number $\alpha$ in $I$ can be written as $\frac{p}{q} + t|I|$ where $t \in [-1, 1]$. Now,

$$|1 - e^{2\pi i q \alpha}| = |1 - e^{2\pi i (q\alpha - p)}| < 4\pi|q\alpha - p| <$$
$$4\pi \left| q \left( \frac{p}{q} + t|I| \right) - p \right| < 4\pi q|I| < \frac{\varepsilon}{MCc}.$$

Hence, $\max_{\alpha \in I} |\tilde{F}_n(\alpha)|_r < \varepsilon/M$. Let us estimate $\max_{\alpha \in I} |\hat{F}_n(\alpha)|_r$:

$$\max_{\alpha \in I} |\Phi^{-1} \circ R_\alpha \circ \Phi - R_\alpha|_r \leq |D\Phi^{-1}|_{r_1} \max_{\alpha \in I} |R_\alpha \circ \Phi - \Phi \circ R_\alpha|_r < M\frac{\varepsilon}{M} = \varepsilon.$$

□

*Proof of Proposition 2.1.* As before, we denote $\Phi_n^{-1} \circ R_\alpha \circ \Phi_n - R_\alpha$ by $\hat{F}_n(\alpha)$, and set

$$r_1 = \max_{A^r} |\mathrm{Im}(T)|, \quad r_2 = r_1 + \max_{\alpha \in [0,1]} |\hat{F}_n(\alpha)|_{r_1}.$$

We can rewrite the definition (2.1) of $G_n(\alpha)$ in the following way:

$$G_n(\alpha) = T_n^{-1} \circ R_\alpha \circ T_n = T^{-1} \circ \Phi_n^{-1} \circ R_\alpha \circ \Phi_n \circ T = T^{-1} \circ (R_\alpha + \hat{F}_n(\alpha)) \circ T,$$

and then
$$|G_n(\alpha) - G(\alpha)|_r = |T^{-1} \circ (R_\alpha + \hat{F}_n) \circ T - T^{-1} \circ R_\alpha \circ T|_r \leq |DT^{-1}|_{r_2}|\hat{F}_n(\alpha)|_{r_1}.$$

Denote the first factor in this product by M (it does not depend on $\alpha$). Then, by Lemma 2.1, we can choose an interval $I$ centered at $\frac{p_n}{q_n}$ so that
$$\max_{\alpha \in I} |\hat{F}_n(\alpha)|_{r_1} < \frac{\varepsilon}{M}.$$

Then $\max_{\alpha \in I} |G_n(\alpha) - G(\alpha)|_r < \varepsilon$.  □

Proof of the corollary is easy, and we leave it to the reader.

*Remark 2.1.* In the proposition above one can estimate the necessary upper bound for the length of the interval $I_n$. However, for our purposes it is sufficient to establish the existence of a "sufficiently small" interval with the required properties, and we do not write out the bounds explicitly. We note only that $\frac{1}{|I_n|}$ in our construction is of superexponential order in $q_n$.

## 3. ANALYTIC NON-LINEARIZABLE DIFFEOMORPHISM

*Proof of Theorem 1.1.* We begin by proving the **convergence of** $G_n(\alpha)$ for an appropriate $\alpha$. Denote $r_n = 10^n, \varepsilon_n = \frac{1}{10^{n+1}}$, and set $I_0 = [0,1]$, $G_0(\alpha) = R_\alpha$. Here we take $c_n > 0$ arbitrary, with the only condition that $\sum c_n = \infty$.

Let us explain the choice of parameters $p_n$ and $q_n$, guaranteeing the existence and non-linearizability of the analytic limit $G = \lim G_n$. Suppose that for $j = 1, \ldots n-1$ the intervals $I_j \subset I_{j-1}$ centered at $\frac{p_j}{q_j}$ are chosen. Pick any $\frac{p_n}{q_n} \in I_{n-1}$ with a large denominator:

$$(3.1) \qquad q_n > 10^{n+2} \sum_{i=1}^{n-1} c_i q_i.$$

To simplify the computations, assume also that

$$(3.2) \qquad q_n = 4s_n q_{n-1}$$

for a large natural $s_n$. We shall choose $I_n$ centered at $\frac{p_n}{q_n}$ in the following way. For $m \leq n$, let $I_{m,n}$ denote an interval, given by Proposition 2.1 with $\varepsilon = \varepsilon_n$, $r = r_n$, $T = T_{m,n-1}$, $G(\alpha) = G_{m,n-1}(\alpha)$ (here $G_{n,n-1}(\alpha) = R_\alpha$), see (1.5).

Let $I_n^1$ denote an interval centered at $\frac{p_n}{q_n}$ given by Corollary 2.1 with $\varepsilon = \varepsilon_n$, $r = 0$, $T = T_{n-1}$, $G(\alpha) = G_{n-1}(\alpha)$ and $\tau = q_n$. Define

$$I_n = \left(\bigcap_{m=1}^{n} I_{m,n}\right) \cap I_n^1.$$

Assume also that the $I_n$ lies strictly inside $I_{n-1}$, and $|I_n| < \frac{1}{q_n^2}$.

Suppose that for all $n$ the parameters are chosen as above and $\hat{\alpha} = \bigcap_n I_n = \lim_n \frac{p_n}{q_n}$. Then we have obtained an infinite number of converging sequences of diffeomorphisms $(G_{m,n}(\hat{\alpha}))_{n=m}^{\infty}$, $\lim_{n \to \infty} G_{m,n}(\hat{\alpha}) = G_{(m)}(\hat{\alpha})$, such that

$$|R_{\hat{\alpha}} - G_{(m)}(\hat{\alpha})|_{r_m} < \frac{1}{10^m}.$$

Indeed, for any fixed $n$, $\hat{\alpha}$ satisfies the conditions of Proposition 2.1 with $\varepsilon = \varepsilon_n$, $r = r_n$, $T = T_{m,n-1}$, $G(\hat{\alpha}) = G_{m,n-1}(\hat{\alpha})$ for any $m \leq n$. By this proposition, for all $m \leq n$, we have
$$|G_{m,n-1}(\hat{\alpha}) - G_{m,n}(\hat{\alpha})|_{r_n} < \varepsilon_n.$$

Fix an arbitrary $m$. For $n \geq m$, the diffeomorphisms $G_{m,n}(\hat{\alpha})$ form a Cauchy sequence in the metric $|\cdot|_{r_m}$, and the limit $G_{(m)}(\hat{\alpha}) \in \overline{\mathcal{A}_{\hat{\alpha}}^{r_m}}$ exists. Moreover,

$$|R_{\hat{\alpha}} - G_{(m)}(\hat{\alpha})|_r \leq \sum_{j=-1}^{\infty} |G_{m,m+j}(\hat{\alpha}) - G_{m,m+j+1}(\hat{\alpha})|_{r_m} <$$
$$< \frac{1}{10^{m+1}} \sum_{j=0}^{\infty} \frac{1}{10^j} < \frac{1}{10^m}.$$

**Non-linearizability.** Now we show that for any $m$ the limit $G_{(m)}(\hat{\alpha})$ of the above construction is non-linearizable. By Remark 1.1, it is enough to show it for $G(\hat{\alpha}) = G_{(1)}(\hat{\alpha})$. Suppose the contrary, i.e. that there exists a homeomorphism $T$ of the torus, such that $G(\hat{\alpha}) = T^{-1} \circ R_{\hat{\alpha}} \circ T$. Then any invariant curve $\Gamma = \Gamma(y_0)$ of $G(\hat{\alpha})$ has the form $T^{-1}\{y = y_0\}$ for some $y_0 \in [0,1]$, and hence is a continuous closed curve, not dense on the torus. We shall come to the contradiction with this assumption by showing that $\Gamma$ is dense on the torus. (In fact, the arguments below imply even more: that the limit diffeomorphism is minimal.) In order to do this, for any fixed $n$ consider invariant curves $\Gamma_n(y_0)$ of $G_n(\hat{\alpha})$. Each of them has the form $T_n^{-1}\{y = y_0\}$, $y_0 \in [0,1]$. Recalling that $T_n^{-1} = \Phi_1^{-1} \circ \ldots \Phi_n^{-1}$, we compute:

$$T_n^{-1}(x,y) = \left(x - \sum_{j=1}^n \frac{1}{2q_j}, y - \sum_{j=1}^n c_j \sin 2\pi q_j (x - \sum_{i=j}^n \frac{1}{2q_i})\right)$$

(3.3)
$$= \left(x - \sum_{j=1}^n \frac{1}{2q_j}, y - \sum_{j=1}^n c_j \sin 2\pi q_j (x - \sum_{i=1}^n \frac{1}{2q_i})\right);$$

the last equality follows from (3.2). Then the lift of $\Gamma_n(y_0)$ is the graph of the function

(3.4)
$$y = \Gamma_n(x) = y_0 - \sum_{j=1}^n c_j \sin 2\pi q_j x.$$

We shall prove the following:
a) for any $\varepsilon$ there exists an $n$ such that the set of points $\{((x + i\hat{\alpha}), \Gamma_n(x + i\hat{\alpha})) \mid i = 0, \ldots q_{n+1}\}$ is $\varepsilon$-dense on $\mathbb{T}^2$;
b) in these points the curve $\Gamma$ is $\varepsilon$-close to $\Gamma_n$, in other words,

$$\max_{i=0,\ldots q_{n+1}} |G^i(\hat{\alpha}) - G_n^i(\hat{\alpha})|_0 < \varepsilon.$$

This will insure that for an arbitrary $\varepsilon$ the curve $\Gamma$ is $\varepsilon$-dense on $\mathbb{T}^2$, contradicting the assumption.

To prove a), fix an $m$ such that $\frac{1}{q_m} \leq \varepsilon$, and consider the intervals $J_k = \left[\frac{k}{q_m}, \frac{k+1}{q_m}\right]$, $k = 0 \ldots q_m - 1$. First we find an $n$ such that $|\max_{J_k} \Gamma_n - \min_{J_k} \Gamma_n| > 1$; this will imply that the curve $\Gamma_n$ is $\varepsilon$-dense on the torus. For any $k \leq q_m - 1$, consider two points in $J_k$: $z_k = \frac{4k+1}{4q_m} + \sum_{i=m+1}^{\infty} \frac{1}{4q_i}$, and $w_k = \frac{4k+3}{4q_m} - \sum_{i=m+1}^{\infty} \frac{1}{4q_i}$. Using (3.1), one verifies that $z_k, w_k \in J_k$. For an arbitrary $n > m$, using (3.2), we have

$$\Gamma_n(z_k) = y_0 - \sum_{j=1}^{n} c_j \sin 2\pi q_j \left(\frac{4k+1}{4q_m} + \sum_{i=m+1}^{\infty} \frac{1}{4q_i}\right) =$$

$$C_1(m) - \sum_{j=m}^{n} c_j \sin\left(\frac{\pi}{2} + 2\pi \sum_{i=j+1}^{\infty} \frac{q_j}{4q_i}\right) = C_1(m) - \sum_{j=m}^{n} c_j \cos\left(2\pi \sum_{i=j+1}^{\infty} \frac{q_j}{4q_i}\right).$$

Here $C_1(m)$ is a constant independent of $n$. By (3.1), $2\pi \sum_{i=j+1}^{\infty} \frac{q_j}{4q_i} \leq \frac{1}{10^j}$. Using the estimate $\cos x \geq 1 - x^2$, valid for small $x$, we have:

$$\Gamma_n(z_k) < C_1(m) - \sum_{j=m}^{n} c_j \left(1 - \frac{1}{10^{2j}}\right) < C_1(m) + \frac{1}{10^m} - \sum_{j=m}^{n} c_j = \tilde{C}_1(m) - \sum_{j=m}^{n} c_j.$$

A similar calculation shows that $\Gamma_n(w_k) > \tilde{C}_2(m) + \sum_{j=m}^{n} c_j$. Since $\sum_{j=1}^{\infty} c_j$ diverges, there exists an $n$ such that $(\Gamma_n(w_k) - \Gamma_n(z_k)) > 1$ for all $k \leq q_m$. Then the curve $\Gamma_n$ is $\varepsilon$-dense on $\mathbb{T}^2$. We take $n$ such that $\frac{1}{10^n} < \varepsilon$.

Since, in particular, $\hat{\alpha} \in I_{n+1}$, then $|\hat{\alpha} - \frac{p_{n+1}}{q_{n+1}}| < \frac{1}{q_{n+1}^2}$ by assumption. This implies that the set of points $\{(x + i\hat{\alpha}) \mid i = 0, \ldots q_{n+1}\}$ is $\frac{2}{q_{n+1}}$-dense on $\mathbb{T}^1$. By (3.1), $\frac{2}{q_{n+1}} < \frac{1}{10^{n+2} \sum_{i=1}^{n} c_i q_i}$; by (3.4),

$$|\Gamma'_n|_0 \leq 2\pi \sum_{i=1}^{n} c_i q_i,$$

and therefore the set $\{((x + i\hat{\alpha}), \Gamma_n(x + i\hat{\alpha})) \mid i = 0, \ldots q_{n+1}\}$ is $\frac{1}{10^n}$-dense on the curve $\Gamma_n$. This proves a).

Condition b) holds since $\hat{\alpha} \in \bigcap_n I_n^1$. Indeed,

$$\max_{i=0,\ldots q_{n+1}} |G^i(\hat{\alpha}) - G_n^i(\hat{\alpha})|_0 \leq \sum_{j=n}^{\infty} \max_{i=0,\ldots q_{n+1}} |G_j^i(\hat{\alpha}) - G_{j+1}^i(\hat{\alpha})|_0 < \frac{1}{10^n} < \varepsilon.$$

Hence, $\Gamma(\hat{\alpha})$ is dense on the torus, contradicting the assumed linearizability of $G(\hat{\alpha})$. $\square$

## 4. Minimal non-ergodic diffeomorphism

*Proof of Theorem 1.2.* To simplify the calculations, suppose that (3.2) still holds for all $n$. Then $T_n(x, y) = \left(x + \sum_{j=1}^{n} \frac{1}{2q_j}, y + \sum_{j=1}^{n} c_j \sin 2\pi q_j x\right)$, $T_n^{-1}(x, y)$ has the form (3.3), and the diffeomorphism $G_n(\alpha)$, given by (1.3), has the form

$$G_n(\alpha; x, y) = (x + \alpha, y + g_n(x) - g_n(x + \alpha)),$$

where $g_n(x) = \sum_{j=1}^n c_j \sin 2\pi q_j x$.

Suppose that $c_n = \frac{1}{n}$ for all $n \in \mathbb{N}$, and an irrational $\alpha = \lim \frac{p_n}{q_n}$ is chosen as in Theorem 1.1. Then the analytic limit $G(\alpha) = \lim G_n(\alpha)$ exists, and there exist diffeomorphisms with the same ergodic properties arbitrary close to $R_\alpha$ in any $r$-metric. Let us study the obtained diffeomorphisms closer.

Since $\{\frac{1}{n}\}_n \in l^2(\mathbb{Z})$, $g_n$ converge almost everywhere to a 1-periodic $L^2$-function $g$ (whose Fourier series is $\tilde{g} = \sum_{j=1}^\infty \frac{1}{j} \sin 2\pi q_j x$). Then, almost everywhere,

$$G(\alpha; x, y) = (x + \alpha, y + g(x) - g(x + \alpha)).$$

It is evident that $G$ is not ergodic with respect to the Lebesgue measure. Indeed, it has an infinite number of ergodic components, supported a.e. on the "strips" of the form

$$\{(x, y - g(x)) \mid x \in \mathbb{T}^1, y \in [y_1, y_2] \subset \mathbb{T}^1\},$$

each having nonzero Lebesgue measure.

A short study of this type of diffeomorphisms can be found in [KH]. In particular, the following statement is proven.

**Lemma 4.1.** *Consider an analytic mapping $f : (x, y) \mapsto (x + \alpha, y + \varphi(x))$, where $\varphi : S^1 \to \mathbb{R}$. Then either $\varphi(x) = g(x) - g(x + \alpha) + r$ for some continuous function $g : S^1 \to \mathbb{R}$ and $r \in \mathbb{R}$, or $f$ is minimal.*

Let us give another, a more elegant proof of the minimality of the limit diffeomorphism (in addition to the one contained in the proof of Theorem 1.1). Assume that $q_n$-s grow faster than a geometric progression (this does not contradict assumptions of Theorem 1.1; in fact, Theorem 1.1 requires much faster growth of $q_n$-s). Then the series $\tilde{g}$ is a so-called *lacunary series*. For such series we have the following general statement [Ka].

**Lemma 4.2.** *If a lacunary series $\sum c_j e^{2\pi i q_j x}$ is a Fourier series of a bounded function, then $\sum |c_j| < \infty$.*

This lemma implies implies that $g$ is unbounded, and therefore, non-continuous.

Suppose that $G(\alpha)$ is not minimal. Then, by Lemma 4.1, for the function $\varphi(x) = g(x) - g(x + \alpha)$ there exists a representation $\varphi(x) = h(x) - h(x + \alpha) + r$ with a continuous function $h = S^1 \to \mathbb{R}$ and $r \in \mathbb{R}$. Integrating the equality $g(x) - g(x + \alpha) = h(x) - h(x + \alpha) + r$, we get that $r = 0$. Now denote $g - h$ by $f$. Then $f(x + \alpha) = f(x)$ for any $x$, implying that $f = const$, by the unique ergodicity of the rotation by $\alpha$. This proves that $g = h + const$, which contradicts the continuity of $h$. Hence, $G(\alpha)$ is minimal.

Then, as remarked in the beginning of the proof, for any positive $r$ and $\varepsilon$ there exists a diffeomorphism $G_{r,\varepsilon}(\alpha)$ such that $|G_{r,\varepsilon}(\alpha) - R_\alpha|_r < \varepsilon$, and $G_{r,\varepsilon}(\alpha) = S^{-1} \circ G(\alpha) \circ S$ with an analytic area preserving diffeomorphism $S$. The latter implies that $G_{r,\varepsilon}(\alpha)$ is minimal and non-ergodic. This finishes the proof of Theorem 1.2. □

## 5. Invariant measures

Here we study the properties of the diffeomorphisms $G_n(\alpha)$, given by formula (1.3), namely, their ergodic invariant measures and the convergence of Birkhoff sums. The results of this section will be used in the proof of Theorem 1.3. Any invariant curve $\Gamma_n$ of $G_n(\alpha)$ is the graph of the function $y = \Gamma_n(x)$ (see (3.4)). For an arbitrary $n$, let us fix an invariant curve $\Gamma_n$. Suppose that $\alpha$ is irrational, and consider the (unique) ergodic invariant probability measure $\mu_n$ of $G_n(\alpha)$, supported on this curve. We begin the section by discussing the dependence of the mean values over $\mu_n$ on the parameters of the construction. Let us denote $\int_{\mathbb{T}^2} f(z) d\mu_n$ by $\hat{f}_n$, and $\int_{\mathbb{T}^2} f(z) dz$ by $\hat{f}$. We prove the following result.

**Lemma 5.1.** *Let $\Gamma_n$ be an invariant curve of $G_n(\alpha)$, $\alpha \in [0,1] \setminus \mathbb{Q}$, and consider the invariant probability measure $\mu_n$ of $G_n(\alpha)$, supported on this curve. Then for any fixed trigonometric polynomial $f$ and $\varepsilon > 0$ there exist numbers $\tilde{c}_n = \tilde{c}_n(\varepsilon, f)$ and $\tilde{q}_n = \tilde{q}_n(\varepsilon, f, c_1, \ldots, c_n, q_1, \ldots q_{n-1})$ such that if $c_n > \tilde{c}_n$ and $q_n > \tilde{q}_n$, then*

$$|\hat{f}_n - \hat{f}| = \left| \int_{\mathbb{T}^2} f(z) d\mu_n - \int_{\mathbb{T}^2} f(z) dz \right| < \varepsilon.$$

*In particular, $\tilde{c}_n$ and $\tilde{q}_n$ are independent of the choice of the invariant curve $\Gamma_n$.*

*Remark 5.1.* Though the diffeomorphism $G_n(\alpha)$ depends on $\alpha$, neither $\mu_n$, nor $\Gamma_n$ depend on it. Hence, $\tilde{c}_n$ and $\tilde{q}_n$ do not depend on $\alpha$.

*Proof of Lemma 5.1.* Fix an arbitrary $n$. It is sufficient to prove the statement for exponential functions, the result for trigonometric polynomials follows by linearity. By the definition of $\mu_n$, for any $k, l \in \mathbb{Z}$, we have:

$$\int_{\mathbb{T}^2} e^{2\pi i(kx+ly)} d\mu_n = \int_0^1 e^{2\pi i(kx+l\Gamma_n(x))} dx.$$

For $k = l = 0$ the latter equals 1, which coincides with the integral of this function over the Lebesgue measure. We shall show that for any non-zero pair $(k, l)$ the integral above is small (and hence close to the Lebesgue measure integral of the corresponding function). Let us consider $k = 0, l = 1$; the general situation can be studied in the same way. For any $\varepsilon > 0$, we shall show that

(5.1) $$\left| \int_0^1 e^{2\pi i \Gamma_n(x)} dx \right| < \varepsilon,$$

provided that $c_n$ and $q_n$ are large enough. Rewrite formula (3.4) for the invariant curve of $G_n(\alpha)$ as

(5.2) $$\Gamma_n(x) = -c_n \sin 2\pi q_n x + \Gamma_{n-1}(x).$$

Denote $|\Gamma'_{n-1}|_0$ by $\beta_n$. Evidently, $\beta_n \leq 2\pi \sum_{j=1}^{n-1} c_j q_j$.

In order to prove the estimate (5.1), we shall for each $n$ approximate $\Gamma_n$ by a piecewise linear curve $y = L_n(x)$ in the following way. Suppose that $c_n \geq 1$ is chosen so large that there exists a natural number $s_n > e^2$ such that

(5.3) $$2^4 \sqrt{c_n/\varepsilon} < s_n, \text{ and } s_n \ln s_n < c_n \varepsilon / 2^4.$$

Suppose that

(5.4) $$q_n > \tilde{q}_n := \beta_n s_n.$$

Divide the $x$-interval $[0,1]$ into $4s_n q_n$ subintervals of equal length $\delta_n := \frac{1}{4s_n q_n}$. We call the intervals $\Delta_j = [(j-1)\delta_n, j\delta_n]$, where $j = 1, \ldots, 4s_n q_n$. The graph of $L_n(x)$ consists of the line segments, connecting the points $((j-1)\delta_n, \Gamma_n((j-1)\delta_n))$ with $(j\delta_n, \Gamma_n(j\delta_n))$.

Let $\alpha_j$ be the inclination of $L_n$ on the interval $\Delta_j$, and let $\gamma_j = \Gamma_n((j-1)\delta_n)$, so that $L_n(x)$ has the form

$$L_n(x) = \alpha_j(x - (j-1)\delta_n) + \gamma_j, \quad x \in \Delta_j.$$

First we show that our approximation is "good enough", i.e. the difference

$$\left| \int_0^1 e^{2\pi i \Gamma_n(x)} dx - \int_0^1 e^{2\pi i L_n(x)} dx \right| =$$

$$= \left| \int_0^1 e^{2\pi i L_n(x)} \left( e^{2\pi i (\Gamma_n(x) - L_n(x))} - 1 \right) dx \right| < 4\pi |\Gamma_n - L_n|_0$$

is less than $\varepsilon/2$. It can be computed (applying the Mean Value theorem on each of the subintervals $\Delta_j$ and using (5.4) and the first inequality of (5.3)) that

$$4\pi |\Gamma_n - L_n|_0 < 4\pi |\Gamma_n''|_0 \delta_n^2 \leq \pi^3 \left( \frac{c_n}{s_n^2} + \frac{\beta_n}{q_n s_n^2} \right) \leq 2\pi^3 \frac{c_n}{s_n^2} < \varepsilon/2.$$

Now we estimate (5.1) with $L_n$ instead of $\Gamma_n$. Note that the inclinations $\alpha_j$ mostly depend on the first term of the sum (5.2). Let us first estimate $\alpha_j$ on the first quarter of the period of $\sin(2\pi q_n x)$, i. e. on $\Delta_j$ for $j = 1, \ldots s_n$. By the Mean Value theorem, for some $z_j \in \Delta_j$, $|\alpha_j| = |\Gamma_n'(z_j)| \geq 2\pi c_n q_n \min_{z \in \Delta_j} |\cos(2\pi q_n z)| - \beta_n$. Then for $j = 1, \ldots (s_n - 1)$ the inclination $\alpha_j$ of $L_n$ on $\Delta_j$ satisfies the following:

$$|\alpha_j| \geq 2\pi c_n q_n \cos(2\pi q_n j \delta_n) - \beta_n = 2\pi c_n q_n \cos \frac{\pi j}{2s_n} - \beta_n$$

$$\geq 2\pi c_n q_n (1 - \frac{j}{s_n}) - \beta_n \geq c_n q_n (1 - \frac{j}{s_n}) := A_j.$$

We used the fact that in this area $\cos x \geq 1 - \frac{2}{\pi} x$, and the estimate $\beta_n \leq c_n q_n (1 - \frac{j}{s_n})$ for all $j = 1, \ldots (s_n - 1)$, which follows from (5.4) and $c_n \geq 1$.

The value of $\alpha_{s_n}$ has to be estimated separately. Using (5.4),

$$|\alpha_{s_n}| = \frac{|\Gamma_n(s_n \delta_n) - \Gamma_n((s_n - 1)\delta_n)|}{\delta_n} \geq \frac{c_n(2\pi q_n \delta_n)^2}{\delta_n} - \beta_n > \frac{c_n q_n}{s_n} := A_{s_n}.$$

It is not difficult to see that if we set

$$A_{2s_n - j + 1} := A_j, \quad j = 1, \ldots s_n,$$

then $|\alpha_j| \geq A_j$ for $j = 1, \ldots 2s_n$.

Now we have got lower bounds on $\alpha_j$ on the interval $[0, \frac{1}{2q_n}]$. Though the curve $L_n$ is not $\frac{1}{2q_n}$-periodic, the estimates of $|\alpha_j|$ can be extended $\frac{1}{2q_n}$-periodically to the whole interval $[0,1]$, since the "non-periodic part" is uniformly bounded by $\beta_n$, and we have the same estimates on each period. We set

$$A_{2s_n k + j} := A_j, \quad j = 1, \ldots 2s_n, \ k = 1, \ldots (2q_n - 1).$$

Then on each $\Delta_j$ we have
$$|\alpha_j| \geq A_j.$$

The preceding estimates can be summarized as follows:
$$\int_0^1 e^{2\pi i L_n(x)} \mathrm{d}x = \sum_{j=1}^{4s_n q_n} \int_{(j-1)\delta_n}^{j\delta_n} e^{2\pi i(\alpha_j(x-(j-1)\delta_n)+\gamma_j)} \mathrm{d}x$$
$$= \sum_{j=1}^{4s_n q_n} e^{2\pi i \gamma_j} \frac{1}{2\pi i \alpha_j} \left(e^{2\pi i \alpha_j \delta_n} - 1\right).$$

The absolute value of the expression above is less than
$$\sum_{j=1}^{4s_n q_n} \frac{1}{|\alpha_j|} = \sum_{j=1}^{4s_n q_n} \frac{1}{A_j} = 4q_n \sum_{j=1}^{s_n} \frac{1}{A_j} \leq 4q_n \left(\frac{s_n}{c_n q_n} + \sum_{j=1}^{s_n-1} \frac{1}{c_n q_n (1 - j/s_n)}\right)$$
$$= 4\frac{s_n}{c_n}\left(1 + \sum_{j=1}^{s_n-1} \frac{1}{s_n - j}\right) = 4\frac{s_n}{c_n}\left(1 + \sum_{j=1}^{s_n-1} \frac{1}{j}\right) \leq 4\frac{s_n}{c_n}(2 + \ln s_n) < \varepsilon/2.$$

The last inequality follows from the second inequality of (5.3). This finishes the proof of the lemma. $\square$

We pass to the study of the convergence of Birkhoff sums of $G_n(\alpha)$ and its dependence on $\alpha$. Recall that the invariant curves of $G_n(\alpha)$ do not depend on $\alpha$. Since for an irrational $\alpha$ the restriction of $G_n(\alpha)$ to any of its continuous invariant curves $\Gamma_n$ is uniquely ergodic, for any fixed $f \in C^0$ the limit

(5.5) $$\lim_{k \to \infty} \frac{1}{k+1} \sum_{i=0}^{k} f \circ G_n^i(\alpha; z) = \int_{\mathbb{T}^2} f(z) \mathrm{d}\mu_n = \hat{f}_n(z)$$

exists for every $z \in \mathbb{T}^2$ (though the limit function is a constant on every invariant curve, it is not necessarily a global constant). The following statement shows that the limit in (5.5) is uniform in $z$ and $\alpha$ for all $\alpha$ in a small interval.

**Lemma 5.2.** *Suppose that $G_n(\alpha)$ is defined by (1.3). Then for any $f \in C^\infty(\mathbb{T}^2)$, $\varepsilon > 0$, and any interval $I$, there exists a subinterval $I'_n \subset I$ and a natural number $\tau$—both depending on $G_n(\alpha)$, $f$, $\varepsilon$— such that for any $\alpha \in I'_n$ and any $k \geq \tau$,*
$$\left|\frac{1}{k+1} \sum_{i=0}^{k} f \circ G_n^i(\alpha) - \hat{f}_n\right|_0 < \varepsilon.$$

*Proof.* Fix an arbitrary $\alpha_0 \in I \setminus \mathbb{Q}$. Then there exists a $K = K(\alpha_0) \in \mathbb{N}$ such that
$$\left|\frac{1}{k+1} \sum_{i=0}^{k} f \circ G_n^i(\alpha_0; z) - \hat{f}_n(z)\right| < \varepsilon/4 \quad \text{for all } k \geq K, \ z \in \mathbb{T}^2.$$

Indeed, for any particular $z$ such a number exists by the unique ergodicity of $G_n(\alpha)$ on any of its invariant curves; the common $K$ for all $z \in \mathbb{T}^2$ exists by compactness of $\mathbb{T}^2$. Now let $I'_n = \{\alpha \mid |\alpha - \alpha_0| < \frac{\varepsilon}{4(K+1)|\mathrm{D}(f \circ T_n^{-1})|_0}\}$, and let $a \in \mathbb{N}$ be any number $a > \frac{4|f|_0}{\varepsilon}$. We set $\tau = a(K+1)$, and show that these $I'_n$ and $\tau$ provide the conclusion of the lemma. Indeed, by the assumptions on $\alpha$ and $a$, for any $z$ and any $\alpha \in I'_n$ we have

$$\left| \sum_{i=0}^{K} f \circ G_n^i(\alpha; z) - (K+1)\hat{f}_n(z) \right| \leq$$

$$\left| \sum_{i=0}^{K} \left( f \circ G_n^i(\alpha; z) - f \circ G_n^i(\alpha_0; z) \right) \right| + \left| \sum_{i=0}^{K} (f \circ G_n^i(\alpha_0; z) - (K+1)\hat{f}_n(z)) \right|$$

$$< (K+1)^2 |\mathrm{D}(f \circ T_n^{-1})|_0 |\alpha - \alpha_0| + (K+1)\frac{\varepsilon}{4} < (K+1)\frac{\varepsilon}{2}.$$

Let $k \geq \tau$, i.e. $k = (K+1)b + j$ with some natural $b \geq a$ and $0 \leq j < (K+1)$. For an arbitrary $z$, denote $G_n^{l(K+1)}(\alpha_0; z)$ by $z_l$ for $l = 0, \ldots b$, and note that $\hat{f}_n(z_l) = \hat{f}_n(z)$ for all $l$. Then we have:

$$\left| \frac{1}{k+1} \sum_{i=0}^{k} f \circ G_n^i(\alpha; z) - \hat{f}_n(z) \right| \leq$$

$$\frac{1}{k+1} \left( \sum_{l=0}^{b-1} \left| \sum_{i=0}^{K} f \circ G_n^i(\alpha; z_l) - (K+1)\hat{f}_n(z) \right| + \left| \sum_{i=0}^{j} f \circ G_n^i(\alpha) \right|_0 + (j+1)|\hat{f}_n|_0 \right)$$

$$\leq \frac{1}{(K+1)b} \left( (K+1)b\frac{\varepsilon}{2} + 2(K+1)|f|_0 \right) < \varepsilon.$$

□

## 6. Unique ergodicity

Let $\mathcal{F} \subset C^\infty(\mathbb{T}^2)$ denote an ordered (countable, dense) set of trigonometric polynomials, and let $\mathcal{F}_n$ be the finite set, consisting of the first $n$ elements of $\mathcal{F}$.

To prove unique ergodicity of $G(\alpha)$ (for an appropriate $\alpha$), it is enough to show (see Th 1.9.2 [M]) that for any function $f \in \mathcal{F}$ and $z \in \mathbb{T}^2$,

$$(6.1) \qquad \lim_{k \to \infty} \frac{1}{k+1} \sum_{i=0}^{k} f \circ G^i(\alpha; z) = \int_{\mathbb{T}^2} f(z) dz = \hat{f}.$$

*Proof of Theorem 1.3.* Here we modify the construction of Theorem 1.1 in such a way that the non-linearizable analytic limit $G(\alpha)$, existing by this theorem, satisfies the above property. Namely, we shall pay more attention to the choice of the growing sequences $c_n$ and $q_n$, and impose even sharper conditions on the decay of the $|I_n|$ than it is required in Theorem 1.1. The proof of Theorem 1.3 occupies the rest of this section, and splits into three steps.

*Step 1.* For every $n$ let us suppose that $c_n \geq \max_{f \in \mathcal{F}_n} \tilde{c}_n(\frac{1}{10^{n+1}}, f)$, and $q_n \geq \max_{f \in \mathcal{F}_n} \tilde{q}_n(\frac{1}{10^{n+1}}, f, c_1, \ldots c_n, q_1, \ldots q_{n-1})$, where $\tilde{q}_n$ and $\tilde{c}_n$ are from Lemma 5.1. Then, by Lemma 5.1, for all $n$ we have

$$(6.2) \qquad |\hat{f}_n - \hat{f}|_0 < \frac{1}{10^{n+1}} \quad \forall f \in \mathcal{F}_n, \; \alpha \in [0, 1].$$

By Lemma 5.2, for any $n$ there exists an interval $I'_n \subset I_n$ and an integer $\tau_n$ such that

$$(6.3) \qquad \left|\frac{1}{k+1}\sum_{i=0}^{k} f \circ G_n^i(\alpha) - \hat{f}_n\right|_0 < \frac{1}{10^{n+1}} \quad \forall\, f \in \mathcal{F}_n,\ \alpha \in I'_n,\ k \geq \tau_n.$$

We use the freedom of choosing $I_{n+1}$ inside $I_n$ in the following way. Let $\frac{p_{n+1}}{q_{n+1}} \in I'_n$. Fix an arbitrary natural $a > 10^{n+1}\max_{f \in \mathcal{F}_{n+1}}|f|_0$. By Corollary 2.1, there exists an interval $I_{n+1} \subset I'_n$, centered at $\frac{p_{n+1}}{q_{n+1}}$, such that for all $\alpha \in I_{n+1}$ and $f \in \mathcal{F}_n$,

$$(6.4) \qquad \max_{i=0,\ldots(\tau_n+1)a}|f \circ G_{n+1}^i(\alpha) - f \circ G_n^i(\alpha)|_0 < \frac{1}{10^{n+1}}.$$

These are the assumptions on $c_n$, $p_n$, $q_n$ and $I_n$ that are sufficient to prove (6.1).

*Step 2.* Assumption (6.4) guarantees that the first $(\tau_n+1)a$ iterates of $f \circ G_{n+1}^i(\alpha)$ are close to those of $f \circ G_n^i(\alpha)$. Let us show that in fact, under the above assumptions, for each $n$ we have

$$(6.5) \qquad \frac{1}{k+1}\left|\sum_{i=0}^{k}(f \circ G_{n+1}^i(\alpha) - f \circ G_n^i(\alpha))\right|_0 < \frac{7}{10^{n+1}} \quad \forall\, f \in \mathcal{F}_n,\ \alpha \in I_{n+1},\ k \in \mathbb{N}.$$

The proof is based on the fact that for large values of $i$ both $f \circ G_n^i(\alpha)$ and $f \circ G_{n+1}^i(\alpha)$ are close to the mean value of $f$.

Let $k = (\tau_n+1)b + j$ with some natural $b \geq a$ and $0 \leq j < (\tau_n+1)$. The left-hand side of (6.5) is less than

$$(6.6) \qquad \left|\frac{1}{k+1}\sum_{i=0}^{k} f \circ G_{n+1}^i(\alpha) - \hat{f}_n\right|_0 + \left|\frac{1}{k+1}\sum_{i=0}^{k} f \circ G_n^i(\alpha) - \hat{f}_n\right|_0.$$

The second term is less than $\frac{1}{10^{n+1}}$ for $k \geq \tau_n$ by (6.3). Let us estimate the first one. Note that for any $z$, we have:

$$\left|\sum_{i=0}^{\tau_n} f \circ G_{n+1}^i(\alpha;z) - (\tau_n+1)\hat{f}_n(z)\right| \leq \left|\sum_{i=0}^{\tau_n} f \circ G_{n+1}^i(\alpha;z) - \sum_{i=0}^{\tau_n} f \circ G_n^i(\alpha;z)\right|$$
$$+ \left|\sum_{i=0}^{\tau_n} f \circ G_n^i(\alpha;z) - (\tau_n+1)\hat{f}_n(z)\right| < \frac{2(\tau_n+1)}{10^{n+1}}.$$

This is true since each of the terms is less than $\frac{\tau_n+1}{10^{n+1}}$ by assumption: the first one by (6.4), the second one by (6.3).

Set $z_l = G_{n+1}^{l(\tau_n+1)}(\alpha;z)$ for $l = 0,\ldots b$. Note that $|\hat{f}_n(z_l) - \hat{f}_n(z)|_0 < 2|\hat{f}_n - \hat{f}|_0 < \frac{2}{10^{n+1}}$ by (6.2) for any $l$. Now the first term of (6.6) is less than

$$\frac{1}{k+1}\left|\sum_{l=0}^{b-1}\sum_{i=0}^{\tau_n} f \circ G_{n+1}^i(\alpha;z_l) + \sum_{i=0}^{j} f \circ G_{n+1}^i(\alpha;z_b) - (k+1)\hat{f}_n(z)\right| \leq$$

$$\frac{1}{k+1}\left(\sum_{l=0}^{b-1}\left|\sum_{i=0}^{\tau_n} f \circ G_{n+1}^i(\alpha;z_l) - (\tau_n+1)\hat{f}_n(z_l)\right| + \left|\sum_{i=0}^{j} f \circ G_{n+1}^i(\alpha)\right|_0 +\right.$$

$$\left.(j+1)|\hat{f}_n|_0\right) + 2|\hat{f}_n - \hat{f}|_0 \leq \frac{1}{b(\tau_n+1)}\left(\frac{2b(\tau_n+1)}{10^{n+1}} + 2\tau_n|f|_0\right) + \frac{2}{10^{n+1}} < \frac{6}{10^{n+1}}.$$

Hence, (6.5) holds true.

*Step 3.* Suppose that the modification, described in Step 1, is done, and $\hat{\alpha} = \bigcap_n I_n \in [0,1] \setminus \mathbb{Q}$. Then, in particular, the analytic non-linearizable limit $G(\hat{\alpha}) = \lim_n G_n(\hat{\alpha})$ exists. We conclude that for any $f \in \mathcal{F}$ and $\varepsilon$ there exists a $\tau = \tau(f, \varepsilon)$ such that for any $k > \tau$

$$\left| \frac{1}{k+1} \sum_{i=0}^{k} f \circ G^i(\hat{\alpha}) - \hat{f} \right|_0 < \varepsilon.$$

Indeed, given arbitrary $f \in \mathcal{F}$ and $\varepsilon > 0$, take $l$ such that $f \in \mathcal{F}_l$ and $\frac{1}{10^l} < \varepsilon$. For an arbitrary $k$ we can write

$$\left| \frac{1}{k+1} \sum_{i=0}^{k} f \circ G^i(\hat{\alpha}) - \hat{f} \right|_0$$
$$\leq \frac{1}{k+1} \left| \sum_{i=0}^{k} (f \circ G^i(\hat{\alpha}) - f \circ G_l^i(\hat{\alpha})) \right|_0 + \left| \frac{1}{k+1} \sum_{i=0}^{k} f \circ G_l^i(\hat{\alpha}) - \hat{f}_l \right|_0$$
$$+ |\hat{f}_l - \hat{f}|_0 = I + II + III.$$

Since $\hat{\alpha} = \bigcap_n I_n$, (6.5) holds for any $n$; hence for any $k$ the first term of this sum estimates as follows:

$$\frac{1}{k+1} \left| \sum_{i=0}^{k} (f \circ G^i - f \circ G_l^i) \right|_0 \leq \sum_{j=l}^{\infty} \frac{1}{k+1} \left| \sum_{i=0}^{k} (f \circ G_j^i - f \circ G_{j+1}^i) \right|_0$$
$$< \frac{7}{10^{l+1}} \sum_{i=0}^{\infty} \frac{1}{10^i} < \frac{8\varepsilon}{10}.$$

Since $\hat{\alpha} \in I_l'$, then, by (6.3), there exists a $\tau_l$ such that for any $k \geq \tau_l$ the term $II$ is less than $\frac{1}{10^{l+1}} < \varepsilon/10$.

Finally, (6.2) guarantees that the term $III$ is less then $\frac{1}{10^{l+1}} < \varepsilon/10$.

Hence, we have proven that (6.1) holds for any function $f \in \mathcal{F}$ and any $z \in \mathbb{T}^2$, which implies the unique ergodicity of $G(\hat{\alpha})$, as it was mentioned in the beginning of the section. Remark 1.1 guarantees the existence, for any $m$, of a uniquely ergodic diffeomorphism $G_{(m)}(\hat{\alpha})$ such that $|G_{(m)}(\hat{\alpha}) - R_\alpha|_{10^m} < \frac{1}{10^m}$. This finishes the proof of Theorem 1.3.

□

## 7. Density and genericity

**Lemma 7.1.** *Suppose that $\alpha = \frac{p_n}{q_n}$, where $|\frac{p_n}{q_n} - \alpha|$ decay sufficiently fast with $n$. Then minimal non-ergodic diffeomorphisms are dense in $\overline{\mathcal{A}_\alpha^r}$ for any $r > 0$, and in $\overline{\mathcal{A}_\alpha}$.*

*The same is true for uniquely ergodic diffeomorphisms.*

*Proof.* It is enough to give a proof only for the first statement, the other one is proved in exactly the same way. Suppose that $|\frac{p_n}{q_n} - \alpha|$ decay sufficiently fast to provide the conclusion of Theorem 1.2. First we show that in this case for any $r > 0$, any $F \in \mathcal{A}_\alpha^r$, and $\varepsilon > 0$ there exists a minimal non-ergodic diffeomorphism $G \in \overline{\mathcal{A}_\alpha^r}$ such that $|F - G|_r < \varepsilon$.

By the definition of $\mathcal{A}_\alpha^r$, for any $F \in \mathcal{A}_\alpha^r$ there exists an analytic measure preserving diffeomorphism $S$ such that $F = S^{-1} \circ R_\alpha \circ S$, and $S^{-1}$ is analytic in some open neighborhood of $\overline{R_\alpha \circ S(A^r)}$. In particular, there exists $\delta$ such that $S^{-1}$ is analytic in $H \circ S(A^r)$ for any analytic diffeomorphism $H$ such that $|(H - R_\alpha) \circ S|_r < \delta$.

Let us denote $r_1 = \sup_{A^r} |\text{Im}(S)|$ (assume that $r_1 > \delta$), and set $M = |DS^{-1}|_{2r_1}$, $\varepsilon_1 = \min\{\delta, \varepsilon/M\}$.

By Theorem 1.2, there exists a minimal non-ergodic analytic diffeomorphism $\tilde{G} \in \overline{\mathcal{A}_\alpha^{r_1}}$ such that $|\tilde{G} - R_\alpha|_{r_1} < \varepsilon_1$. Then $S^{-1}$ is well defined on $\tilde{G} \circ S(A^r)$, and we can set $G = S^{-1} \circ \tilde{G} \circ S$.

Evidently, $G \in \overline{\mathcal{A}_\alpha^r}$, and the desired estimate holds true:

$$|F - G|_r = |S^{-1} \circ R_\alpha \circ S - S^{-1} \circ \tilde{G} \circ S|_r \leq M|\tilde{G} - R_\alpha|_{r_1} < \varepsilon.$$

Hence, for an arbitrary $r$, minimal non-ergodic diffeomorphisms are dense in $\overline{\mathcal{A}_\alpha^r}$. The same statement for $\overline{\mathcal{A}_\alpha}$ follows from the definition of the topology $\tau$ (see Section 1.1). □

**Lemma 7.2.** *The set of real numbers $\alpha$, satisfying the conditions of Theorem 1.1, is generic in $\mathbb{R}$. The same is true for Theorems 1.2 and 1.3.*

*Proof.* Note that the only condition on $\alpha = \frac{p_n}{q_n}$, sufficient to provide the desired theorems, is the "fast" decay of $|\frac{p_n}{q_n} - \alpha|$. It is an easy exercise to show that such numbers $\alpha$ are generic on $[0, 1]$. □

Theorems A and B follow from Lemmas 7.1, 7.2.

At the end we would like to discuss genericity of the unique ergodicity. The following statement is proven in [FH].

**Lemma 7.3.** *In the space of homeomorphisms of a compact metric space, uniquely ergodic homeomorphisms constitute a $G_\delta$-set.*

This implies that for any $r > 0$, in $\overline{\mathcal{A}_\alpha^r}$ endowed with the topology $\tau_r$, the set of uniquely ergodic diffeomorphisms is a $G_\delta$-set. Combined with Lemma 7.1, this provides genericity of uniquely ergodic diffeomorphisms in any $\overline{\mathcal{A}_\alpha^r}$.

*Acknowledgements.* I am very grateful to my research advisor L. H. Eliasson for his continuous help and attention.

Department of Mathematics, Royal Institute of Technology, S-10044, Sweden
*Current address*: Université Paris 7
*E-mail address*: `masha@math.jussieu.fr`